\documentclass[a4paper]{amsart}  

\usepackage{amssymb} \usepackage{amscd} \usepackage{times}

\def\zbb{\mathbb{Z}}  
  
  \def\phi{\varphi}
 \def\p1{{\mathbb{P}^1_\zbb}}

\newtheorem{Theorem}{\quad Theorem}[section]

\newtheorem{Lemma}[Theorem]{\quad Lemma}
\newtheorem{Proposition}[Theorem]{\quad Proposition}

\newcommand{\be} {\begin{equation}}
\newcommand{\ee} {\end{equation}}

\begin{document}

\title{ An estimate for an elliptic equation in dimension 4.}

\author{Samy Skander Bahoura}

\address{Department of Mathematics, Pierre et Marie Curie University, 75005, Paris, France.}
              
\email{samybahoura@gmail.com, samybahoura@yahoo.fr} 

\date{}

\maketitle

\begin{abstract}

We give a uniform estimate for solutions of prescribed scalar curvature type equation in dimension 4.

\end{abstract}

\section{Introduction and Main Result}

We are on a Riemannian manifold $ (M,g) $ of dimension $ n = 4 $. In this paper we denote $ \Delta_g = -\nabla^j(\nabla_j) $ the Laplace-Beltrami operator and $ N=\frac{2n}{n-2}=4 $ the critical Sobolev exponent in dimension 4.

We consider the following equation (of type prescribed scalar curvature)

\be \Delta_g u+hu=V u^3,\,\, u >0 . \ee

Where $ V $ is a function and $ h $ is a smooth bounded function such that $ h\leq \frac{S_g}{6}$ with $ S_g $ the scalar curvature of $ (M,g) $ and $h$ is bounded in $ C^1 $ norm by $ h_0 \in {\mathbb R}^+ $, a number and $ h\not \equiv \frac{S_g}{6} $.

With the previous conditions, this equation arise in Physics and in Astronomy. Equations of this type were studied by many authors, see, [1-23]. For $ a, b, A >0 $, we consider $ u>0 $ solution of the previous equation relative to $ V $(a Lipschitz function) with the following conditions:

$$ 0 <a \leq V \leq b < +\infty, \,\, {\rm and}\,\, ||\nabla V||_{\infty} \leq A. $$

Our main result is:

\begin{Theorem} For all $ a, b, m >0 $, $ A\to 0 $ and all compact $ K $ of $ M $, there is a positive constant $ c=c(a,b,A,h_0, m, K, M, g) $ such that:

$$ \sup_K u \leq c \,\, {\rm if} \,\, \inf_M u\geq m >0. $$

\end{Theorem}

A consequence of this theorem is: if we take two sequences $ (u_i,V_i) $ with $ u_i $ solution of $(1) $ relative to $ V_i $ with the previous conditions for $a,b>0$ and $ A=A_i\to 0$, we have the following implicit Harnack inequality:

$$ \sup_K u_i\leq c(a,b, h_0, (A_i)_i, \inf_M u_i,K, M, g), $$

and, the function $ c $ change if we consider other sequences $ (v_i,W_i) $ with the same hypothesis for $ W_i $ (Lipschitz and the new Lipschitz constant $ B_i \to 0 $).

\section{Proof of the theorem.}

Let us consider $ x_0 \in M $, by a conformal change of the metric $ \tilde g = \phi^{4/(n-2)} g=\phi^2 g  $ with $ \phi >0 $ we can consider the equation:

\be \Delta_{\tilde g}u + R_{\tilde g} u=Vu^3+(R_g-h)\phi^{-2} u,\,\, u >0 . \ee 

with,

$$ Ricci_{\tilde g}(x_0)= 0. $$

Here; $ R_g= \dfrac{1}{6} S_g $ and  $ R_{\tilde g}= \dfrac{1}{6} S_{\tilde g} $, here we have $ R_g-h\geq 0 $, this function is smooth and uniformly bounded in $ C^1 $ norm.

See the computations in previous papers [3,8], Also, we use the notations of the papers [5,7,8].

\smallskip

1) {\it Blow-up  analysis:} we argue by contradiction and we suppose that $ \sup $ is not bounded. Let $ x_0 \in M $.

\smallskip

We assume that for $ a, b, m >0$ and $ A\to 0$:

\bigskip

$ \forall \,\, c,R >0 \,\, \exists \,\, u_{c,R} $ solution of $ (2) $ such that:

$$ R^{n-2} \sup_{B(x_0,R)} u_{c,R} \geq c\,\, {\rm and} \,\, u_{c,R}\geq m >0, \qquad (H) $$

\bigskip

\begin{Proposition}

There exist a sequence of points $ (y_i)_i $, $ y_i \to x_0 $ and two sequences of positive real number $ (l_i)_i, (L_i)_i $, $ l_i \to 0 $, $ L_i \to +\infty $, such that if we consider $ v_i(y)=\dfrac{u_i[\exp_{y_i}(y)]}{u_i(y_i)} $, we have:

$$ i) \qquad 0 < v_i(y) \leq  \beta_i \leq 2, \,\, \beta_i \to 1. $$

$$ ii) \qquad v_i(y)  \to \dfrac{1}{1+{|y|^2}}, \,\, {\rm uniformly \,\, on\,\, every \,\, compact \,\, set \,\, of } \,\, {\mathbb R}^n . $$

$$ iii) \qquad l_i [u_i(y_i)]  \to +\infty $$
\end{Proposition}

Proof: see [3,8].

\bigskip

2) {\it Polar coordinates and "moving-plane" method}

\smallskip

Let, 

$$ w_i(t,\theta)=e^{t}\bar u_i(e^t,\theta) = e^{t}u_io\exp_{y_i}(e^t\theta), \,\, {\rm and} \,\, a(y_i,t,\theta)=\log J(y_i,e^t,\theta). $$ 

\smallskip

\begin{Lemma}

The function $ w_i $ is solution of:

\be  -\partial_{tt} w_i-\partial_t a \partial_t w_i+\Delta_{\theta}w_i+c w_i=V_iw_i^3+(R_g-h)\phi^{-2} e^{2t}w_i, \ee
 
with,

 $$ c = c(y_i,t,\theta)=1+\partial_t a+R_{\tilde g} e^{2t}. $$ 

\end{Lemma}

Proof of the lemma, see [5,7,8].

\bigskip

Now we have, $ \partial_t a=\dfrac{ \partial_t b_1}{b_1} $, $ b_1(y_i,t,\theta)=J(y_i,e^t,\theta)>0 $,

\bigskip

We can write,

$$ -\dfrac{1}{\sqrt {b_1}}\partial_{tt} (\sqrt { b_1} w_i)+\Delta_{\theta}w_i+[c(t)+ b_1^{-1/2} b_2(t,\theta)]w_i=V_iw_i^3+ (R_g-h)\phi^{-2}e^{2t}w_i, $$

where, $ b_2(t,\theta)=\partial_{tt} (\sqrt {b_1})=\dfrac{1}{2 \sqrt { b_1}}\partial_{tt}b_1-\dfrac{1}{4(b_1)^{3/2}}(\partial_t b_1)^2 .$

\bigskip

Let,

$$ \tilde w_i=\sqrt {b_1} w_i, $$

we have:

\begin{Lemma}

The function $ \tilde w_i $ is solution of:

$$ -\partial_{tt} \tilde w_i+\Delta_{\theta}(\tilde w_i)+2\nabla_{\theta}(\tilde  w_i) .\nabla_{\theta} \log (\sqrt {b_1})+(c+b_1^{-1/2} b_2-c_2) \tilde w_i= $$

\be = V_i\left (\dfrac{1}{b_1} \right ) {\tilde w_i}^3+(R_g-h)\phi^{-2} e^{2t}{\tilde w_i}, \ee

where, $ c_2 =[\dfrac{1}{\sqrt {b_1}} \Delta_{\theta}(\sqrt{b_1}) + |\nabla_{\theta} \log (\sqrt {b_1})|^2].$
\end{Lemma}
Proof of the lemma, see [5,7,8].
\smallskip

We have,

$$ c(y_i,t,\theta)=1+\partial_t a + R_{\tilde g} e^{2t}, \qquad (\alpha_1) $$ 

$$ b_2(t,\theta)=\partial_{tt} (\sqrt {b_1})=\dfrac{1}{2 \sqrt { b_1}}\partial_{tt}b_1-\dfrac{1}{4(b_1)^{3/2}}(\partial_t b_1)^2 ,\qquad (\alpha_2) $$ 

$$ c_2=[\dfrac{1}{\sqrt {b_1}} \Delta_{\theta}(\sqrt{b_1}) + |\nabla_{\theta} \log (\sqrt {b_1})|^2], \qquad (\alpha_3) $$

Then,

$$ \partial_{t}c(y_i,t,\theta)=\partial_{tt}a, $$

by proposition 2.1,

$$ |\partial_tc_2|+|\partial_t b_1|+|\partial_t b_2|+|\partial_t c|\leq K_1e^{2t}. $$

We have for $ \lambda_i= - \log u_i(y_i) $,

$$ w_i(2\xi_i-t,\theta)=w_i[(\xi_i-t+\xi_i-\lambda_i-2)+(\lambda_i+2)] , $$

Thus,
 $$ w_i(2\xi_i-t,\theta)=e^{[2(\xi_i-t+\xi_i-\lambda_i-2)]/2}e^{2}v_i[\theta e^2e^{(\xi_i-t)+(\xi_i-\lambda_i-2)}]
\leq 2e^2=\bar c. $$

\bigskip

3){\it The "moving-plane" method:}

\bigskip

Let $ \xi_i $ a real number,  and suppose $ \xi_i \leq t $. We set $ t^{\xi_i}=2\xi_i-t $ and $ \tilde w_i^{\xi_i}(t,\theta)=\tilde w_i(t^{\xi_i},\theta) $.

\bigskip

We have, 

$$ -\partial_{tt} \tilde w_i^{\xi_i}+\Delta_{\theta}(\tilde w_i)+2\nabla_{\theta}(\tilde  w_i^{\xi_i}) .\nabla_{\theta} \log (\sqrt {b_1}) \tilde w_i^{\xi_i}+[c(t^{\xi_i})+b_1^{-1/2}(t^{\xi_i},.)b_2(t^{\xi_i})-c_2^{\xi_i}] \tilde w_i^{\xi_i}= $$

$$ =V_i\left (\dfrac{1}{b_1^{\xi_i}} \right ) { ({\tilde w_i}^{\xi_i}) }^3+ $$

$$ + (((R_g-h)\phi^{-2})e^{2t})^{\xi_i}{\tilde w_i}^{\xi_i}, $$

By using the same arguments than in [3, 5, 7, 8], we have:

\bigskip

\begin{Proposition}

We have for $ \lambda_i= - \log u_i(y_i) $;

$$ 1)\,\,\, \tilde w_i(\lambda_i,\theta)-\tilde w_i(\lambda_i+4,\theta) \geq \tilde k>0, \,\, \forall \,\, \theta \in {\mathbb S}_{3}. $$

For all $ \beta >0 $, there exist $ c_{\beta} >0 $ such that:

$$ 2) \,\,\, \dfrac{1}{c_{\beta}} e^{(n-2)t/2}\leq \tilde w_i(\lambda_i+t,\theta) \leq c_{\beta}e^{(n-2)t/2}, \,\, \forall \,\, t\leq \beta, \,\, \forall \,\, \theta \in {\mathbb S}_{3}. $$

\end{Proposition}

We set,

$$ \bar Z_i=-\partial_{tt} (...)+\Delta_{\theta}(...)+2\nabla_{\theta}(...) .\nabla_{\theta} \log (\sqrt {b_1})+(c+b_1^{-1/2} b_2-c_2)(...) $$

{\bf Remark:} In the operator $ \bar Z_i $, by using the proposition 3, the coeficient $ c+b_1^{-1/2}b_2-c_2 $ satisfies:

$$ c+b_1^{-1/2}b_2-c_2 \geq k'>0,\,\, {\rm pour }\,\, t<<0, $$

it is fundamental if we want to apply Hopf maximum principle.

\bigskip

\underbar {\it Goal:}

\bigskip

Set $ \bar w_i= \tilde w_i-\frac{m}{2} e^t $. Like in [8] we have the some properties for $\bar w_i $, we have:

\smallskip

\begin{Lemma}

\bigskip

There is $ \nu <0 $ such that for $ \lambda \leq \nu $ :

$$ \bar w_i^{\lambda}(t,\theta)-\bar w_i(t,\theta) \leq 0, \,\, \forall \,\, (t,\theta) \in [\lambda,t_i] \times {\mathbb S}_{3}. $$

\end{Lemma}

\smallskip

Let $ \xi_i $ be the following real number,

$$ \xi_i=\sup \{ \lambda \leq \lambda_i+2, \bar w_i^{\xi_i}(t,\theta)-\bar w_i(t,\theta) \leq 0, \,\, \forall \,\, (t,\theta)\in [\xi_i,t_i]\times {\mathbb S}_{3} \}. $$

\smallskip

Like in [3,5,7,8], we have elliptic second order operator. Here it is $ \bar Z_i $, the goal is to use the "moving-plane" method to have a contradiction. For this, we must have:

$$ \bar Z_i(\bar w_i^{\xi_i}-\bar w_i) \leq 0, \,\, {\rm if} \,\, \bar w_i^{\xi_i}-\bar w_i \leq 0. $$

Clearly, we have:

\begin{Lemma}

$$ b_1(y_i,t,\theta)=1-\dfrac{1}{3} \tilde Ricci_{y_i}(\theta,\theta)e^{2t}+\ldots, $$

$$ R_{\tilde g}(e^t\theta)=R_{\tilde g}(y_i) + <\nabla R_{\tilde g}(y_i)|\theta > e^t+\dots . $$
\end{Lemma}

We have:

$$ \partial_t[(R_g-h)\phi^{-2} e^{2t}] \geq -\epsilon_i e^{2t}, \epsilon_i \to 0.$$

Thus,

$$[(R_g-h)\phi^{-2} e^{2t}]-[(R_g-h)\phi^{-2} e^{2t}]^{\xi_i}\geq -\epsilon_i(e^{2t}-e^{2t^{\xi_i}}) $$

Now, we write:

$$ \bar Z_i(\bar w_i^{\xi_i}-\bar w_i) \leq {\tilde A_i}(e^{t}-e^{t^{\xi}})(\tilde w_i^{\xi_i})^3+V_i(b_1^{\xi_i})^{-1}(\tilde w_i^3-(\tilde w_i^{\xi_i})^3)+o(1)e^{2t}(e^t-e^{t^{\xi_i}})+o(1)\tilde w_i^{\xi_i}(e^{2t}-e^{2t^{\xi_i}}) $$
We can write,

$$ e^{2t}-e^{2t^{\xi_i}}=(e^t-e^{t^{\xi_i}})(e^t+e^{t^{\xi_i}}) \leq 2 e^t (e^t-e^{t^{\xi_i}}). $$

Thus,

$$ \bar Z_i(\bar w_i^{\xi_i}-\bar w_i) \leq 4 e \tilde A_i (e^t-e^{t^{\xi_i}}) ({\tilde w_i}^{\xi_i})^2 + V_i (b_1^{\xi_i})^{-1/2}[( \tilde w_i^{\xi_i})^3-{\tilde w_i}^3]+ o(1) e^{2t} (e^t-e^{t^{\xi_i}})+o(1)e^t {\tilde w_i}^{\xi_i}(e^t-e^{t^{\xi_i}}). $$

But,

$$  0 < {\tilde w_i}^{\xi_i} \leq 2e, \,\,\, \tilde w_i \geq \dfrac{m}{2} e^t \,\,\, {\rm and} \,\,\, {\tilde w_i}^{\xi_i}-{\tilde w_i} \leq \dfrac{m}{2} ( e^{t^{\xi_i}}-e^t), $$

and,

$$ ({\tilde w_i}^{\xi_i})^3-{\tilde w_i}^3=({\tilde w_i}^{\xi_i}-\tilde w_i)[ ({\tilde w_i}^{\xi_i})^2+{\tilde w_i}^{\xi_i} {\tilde w_i}+ {\tilde w_i}^2] \leq ({\tilde w_i}^{\xi_i}-{\tilde w_i}) ({\tilde w_i^{\xi_i}})^2 + ({\tilde w_i}^{\xi_i}-{\tilde w_i}) \dfrac{m^2 e^{2t}}{4} + ({\tilde w_i}^{\xi_i}-{\tilde w_i}) \dfrac{m}{2} e^t {\tilde w_i}^{\xi_i}, $$

then,

$$ \bar Z_i(\bar w_i^{\xi_i}-\bar w_i) \leq    \left [ ({\tilde w_i}^{\xi_i})^2 [\dfrac{a m}{4} - 4 e A_i] + [\dfrac{am^3}{16}-o(1)]e^{2t} + [\dfrac{am^2}{8}- o(1)] e^t {\tilde w_i}^{\xi_i} \right ] ( e^{t^{\xi_i}}- e^t) \leq 0. $$

I fwe use the Hopf maximum principle, we obtain (like in [3,5,7,8]):

$$ \min_{\theta \in {\mathbb S}_3} w_i(t_i,\theta) \leq \max_{\theta \in {\mathbb S}_3} w_i(2\xi_i-t_i, \theta), $$

we can write (by using the proposition 2.4):

$$  l_i u_i(y_i) \leq c, $$

Contradiction.

\end{document}